\documentclass[centertags,leqno]{article}

\usepackage{amsmath}
\usepackage{amssymb}
\usepackage{latexsym}
\usepackage{amsxtra}
\usepackage{amscd}
\usepackage{theorem}
\usepackage{xypic}
\usepackage{graphicx}
\usepackage{mathbbol}
\usepackage{graphicx}
\usepackage{wrapfig}
\usepackage{amscd, wasysym}
\usepackage[matrix,arrow,curve]{xy}

\theoremstyle{plain}
\newtheorem{lemma}{Lemma}[section]
\newtheorem{theorem}[lemma]{Theorem}
\newtheorem{utv}[lemma]{Proposition}
\newtheorem{sled}[lemma]{Corollary}
\theoremstyle{definition}
\newtheorem{definition}[lemma]{Definition}
\newtheorem{zamech}[lemma]{Remark}
\newtheorem{notation}[lemma]{Notation}

\newtheorem{example}[lemma]{Example}

\title{Abel pairs and modular curves}
\author{Dmitry Oganesyan}
\date{}
\begin{document}
\maketitle

\begin{abstract}
In this article we consider rational functions on algebraic curves, which have one zero and one pole (and call pair of such function and curve Abel pair). We investigate moduli spaces of such functions on curves of genus one; the number of Belyi pairs among them is calculated for fields $\mathbb C$ and $\overline{\mathbb F_p}$. This result could be fruitfully used for investigation of Hurwitz's space and modular curves for fields of finite characteristic.
\end{abstract}

\section{Introduction}

In this paper we consider algebraic curves and rational function on them with the divisor of a certain combinatorial type. More specifically, we study pairs $(\mathcal X,\alpha)$, where $\mathcal X$ is an algebraic curve and $\alpha$ is a function on it with $\mathrm{div}(\alpha)=nA-nC$. We call such a pair an \emph{Abel pair}. An Abel function has two critical values $0$ and $\infty$ of valency $n$; if it has only one additional critical value it is a \emph{Belyi function} aslo. One of the main results of this paper is the calculation of the number of Abel-Belyi pairs of a given degree of genus 1 over $\mathbb C$ (Theorem 5.8) and $\overline{\mathbb F_p}$ (Theorem 6.3).

The idea of calculation is to include Abel-Belyi pairs into families of Abel pairs. In this families the Abel-Belyi pairs correspond to zeros of some function on the base of the family. Knowing multiplicities of these zeros and the degree of this function we count the Abel-Belyi pairs of genus 1. This method works in positive characteristic and for other types of Belyi pairs as well. The families used in this construction are interesting by themselves; the bases of the families are modular curves $X_1(n)$.

Similar consideration of such Belyi pairs can be found in ~\cite{Pakovich}.

I express my gratitude to my supervisor, George B. Shabat, and to the  participants of Mechmath MSU seminar ``Graphs on surfaces and curves over number fields'' for their help and fruitful discussion in performing this work. I also thank my reviewer for his advises and corrections.

\section{Abel pairs}
\begin{definition}
An \it Abel pair \rm is a pair $(\mathcal X,\alpha)$, where $\mathcal X$ is a complete smooth algebraic curve over an algebraically closed field $\Bbbk$ and $\alpha$ a rational function on it, whose divisor has the form $\mathrm{div}(\alpha)=nA-nB$.
Such an $\alpha$ will be called an \emph{Abel function}.
\end{definition}

\begin{example}\label{exampleone}

Consider a family of elliptic curves $y^2=(1+kx)^2-4x^3$ ($j$-invariant $j=-{\frac {{k}^{3} \left( {k}^{3}+24 \right) ^{3}}{27+{k}^{3}}}$), with $k\in \mathbb \Bbbk\setminus \{{-3\sqrt[3]{1}}\}$.
It is easy to check that $\alpha=1+kx-y$ is an Abel function on it.
\end{example}

How is this definition related with Abel? Abel described in \cite{Abel} the type of elliptic integrals of the third kind which can be expressed in terms of elementary functions. These integrals are related to Abel pairs. For example~\ref{exampleone} we obtain ($v=1/x$): $$\int\frac{(3u+k)dx}{\sqrt{(u^2+ku)^2-4u}}=\ln(u(u+k)^2-2+(u+k)\sqrt{(ku+u^2)^2-4u})+C $$
This integration is explained by the fact that the difference of two infinite points on the curve $v^2=(u^2+ku)^2-4u$ has order 3.

\begin{definition} Let $ f: \mathcal X_1 \to \mathcal X_2 $ be a morphism of curves and $\lambda\in \operatorname{PSL}_2(\Bbbk)$. We call a pair $\mathcal F=(f,\lambda)$ a \emph{morphism} from an Abel pair $(\mathcal X_1,\alpha_1)$ to an Abel pair $(\mathcal X_2,\alpha_2)$  if $\lambda\circ \alpha_1=\alpha_2\circ f$ i.e. if the diagram
$$
\begin{CD}
X_1 @>f>> X_2 \\
@V{\alpha_1}VV @V{\alpha_2}VV \\
\mathbb P^1(\Bbbk) @>{\lambda}>> \mathbb P^1(\mathbb{\Bbbk})
\end{CD}
$$
is commutative.
\end{definition}

\begin{zamech}
So, pairs $(\mathcal X, \alpha)$, $(\mathcal X, 5\alpha)$ (for $\operatorname{char}\Bbbk\ne 5$) and $(\mathcal X, 1/\alpha)$ are isomorphic.
\end{zamech}

\begin{definition}
We call an Abel pair $(\mathcal X, \alpha)$ an \emph{Abel-Belyi pair} if it is also a Belyi pair, i.e. if $\alpha$ has only three critical values.
\end{definition}

The dessin d'enfant $\beta^{-1}{([0,1])}$ on $\mathcal X$ is associated to a Belyi pair $(\mathcal X, \beta)$ over the field $\mathbb C$. We color vertices $\beta^{-1}(0)$ in black and vertices $\beta^{-1}(1)$ in white, see ~\cite{LandoZvonkin}.

The dessin corresponding to an Abel-Belyi pair has one black vertex and one face.

\section{Embedded graph associated with Abel pair}
\begin{wrapfigure}{r}{1.8in}
\includegraphics[width=1.6in,height=1.9in]{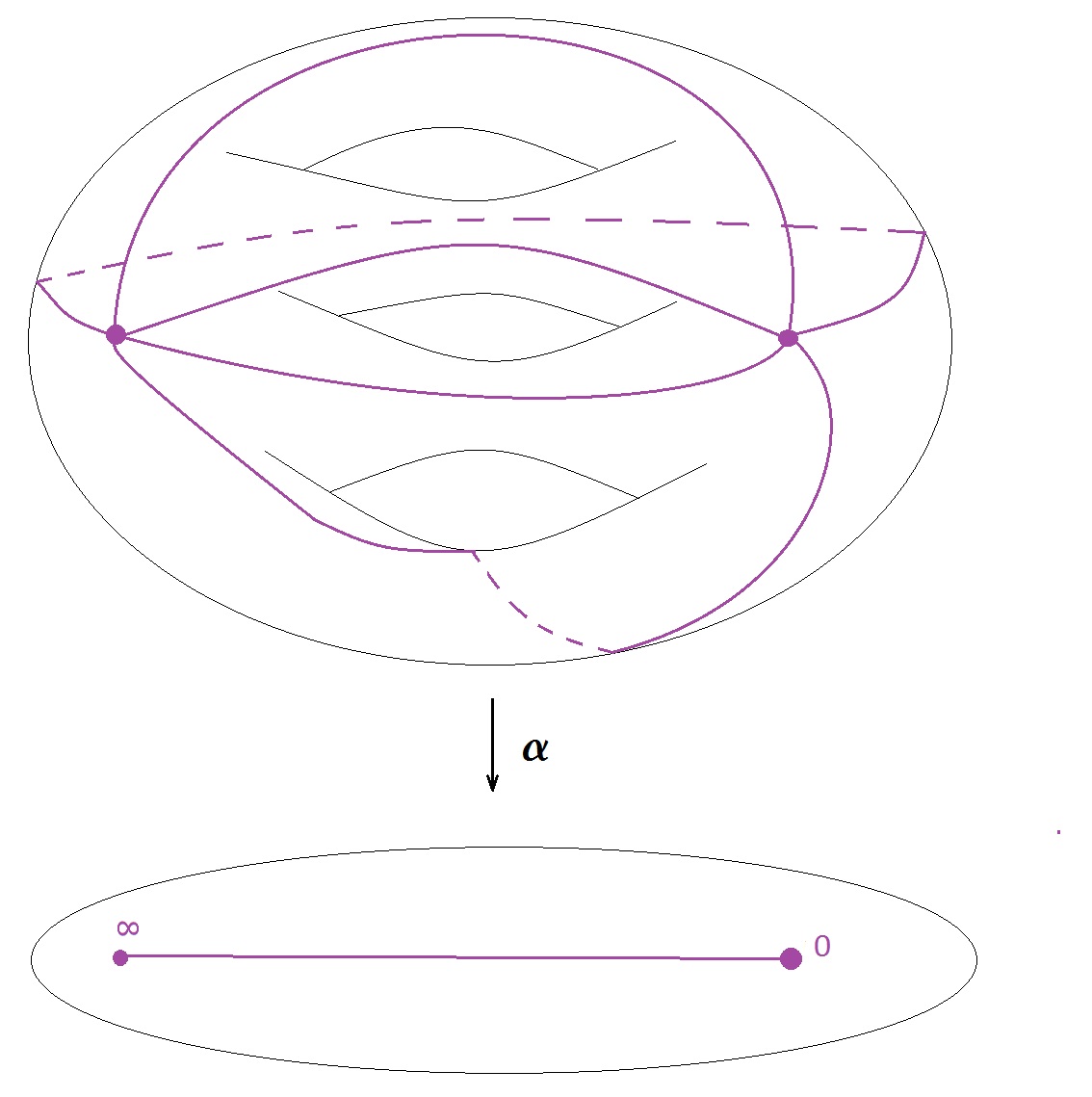}
\end{wrapfigure}
Let $(\mathcal X,\alpha)$ be a generic Abel pair. \emph{Generic} means that $\alpha$ has no critical values in $\mathbb{R}_{<0}$. Denote $$\Gamma_{\mathcal X, \alpha}: = \alpha^{-1}(\mathbb{R}_{\le 0} \cup \{\infty \}) \subset \mathcal X$$

\begin{utv} If $ \alpha $ has no critical values on the half-line $\mathbb {R}_{<0} $, then
$ \Gamma_{\mathcal X,\alpha} $ is embedded in $\mathcal X$; this graph has two vertices and $n =\deg\alpha $ edges.
\end{utv}

$ \Gamma_ {\mathcal X, \alpha} $ is not necessarily a dessin because some of its faces can be not simply-connected. For example, for the pair $ (\mathcal X, \alpha) $, where $\mathcal X$ is given by the equation $y^2=x^3+x$ and $\alpha = x$, the complement $\Gamma_ {\mathcal X, \alpha} $ is not a dessin since $\mathcal X \backslash \Gamma_ {\mathcal X, \alpha} $ is homeomorphic to a cylinder.

Next we consider combinatorics of these embedded graphs. The semi-edges incident to the vertex $\alpha^{-1}(0)$ are cyclically ordered, the same is true for $\alpha^{-1}(\infty)$. The edges of this graph define a bijection between two cyclically ordered $n$-elementary sets of semi-edges. By counting the number of such bijections we get the proposition. (Abel-Belyi pairs are counted with weight $\frac 1{\# Aut}$)

\begin{utv}
The number of Abel-Belyi pairs of degree $n$ is $$m(n)=\sum\limits_{\deg\alpha=n}\frac1{\# Aut(\mathcal X,\alpha)}=\frac{(n-1)!}{n}$$.
\end{utv}

At the proposition above we have counted together all the Abel-Belyi pairs of fixed degree of all possible genera. Next we want to separate them by the genus and the list of valencies. Also we want to establish the similar formulas in positive characteristic. So we undertake some algebraic considerations.

\section{Defining equations of Abel curves}

Let $\mathcal X$ be an irreducible algebraic curve over $\Bbbk$ and $\alpha\in\Bbbk [\mathcal X]$ a function on it, $\deg \alpha=n$, and $\operatorname{char}\Bbbk\nmid n$. Note that $\mathcal X$ may be singular.

\begin{theorem} Suppose that $\alpha:\mathcal X\to \mathbb P_1(\Bbbk)$ is a separable covering. Then $(\mathcal X,\alpha)$ is an Abel pair of degree $n$ if and only if there exist $x\in \Bbbk [\mathcal X]$, $\deg x=k$, and the polynomials $P_1$, $P_2$,$\ldots,P_{k-1}\in \Bbbk[x]$, $ \deg P_i \le n $, such that the equation of the curve $\mathcal X $ can be written as
$$\alpha^k+\alpha^{k-1} P_{k-1}(x)+\cdots+\alpha P_1(x)+x^n=0$$
\end{theorem}

\textbf{Proof:} Let $ (\mathcal X, \alpha) $ be an Abel pair, then $\Bbbk (\mathcal X) / \Bbbk (\alpha) $ is a finite extension of fields. Primitive element theorem says there is an element $x$ such that $\Bbbk(\mathcal X)=\Bbbk(\alpha)(x)$. Let $ F \in \Bbbk [\alpha, t] $ be minimal polynomial of $x$, i.e. $\Bbbk(\mathcal X)=\Bbbk(\alpha)[t]/(F)$.

Let $\operatorname{div} (\alpha) = nA-nB $. Using a linear fractional transformation, we assume that $ x (A) = 0 $ and $x(B) = \infty $.

So we obtain the equation of the curve $ F (\alpha, x)=0 $.

We expand $F$ in powers of $ \alpha $: $\alpha^k P_k (x) + \alpha ^ {k-1} P_{k-1} (x) + \cdots + \alpha
P_1 (x) + P_0 (x) $. Considering the degree of $ \alpha $, we get $ \operatorname{deg} P_i \le n $. The condition $ \operatorname{div} (\alpha) = nA-nB $ implies $P_0 (x) = x^n$ and $P_k (x) = 1$.

Conversely, from the form of the equation it is clear that $ \alpha $ has the divisor $ nA-nB $. $\Box$

Further in this article we will assume that $\mathcal X$ is a hyperelliptic curve. In this case we obtain the following proposition.

\begin{utv} Suppose $\operatorname{char}\Bbbk\ne 2$.
\begin{itemize}
\item[\emph{(i)}] $(\mathcal X,\alpha)$ is a hyperelliptic Abel pair of degree $n$ if and only if there exist $x\in \mathbb C[\mathcal X]$, $\deg x=2$, and polynomial $P$, $ \deg P \le n $, such that the equation of the curve $\mathcal X $ can be written as
$$F(\alpha,x)=\alpha^2+\alpha P(x)+x^n=0$$

\item[\emph{(ii)}] If genus of $X$ is 1 than $\deg P\le \frac n2$, and also in the case of even n the leading coefficient of $P$ is 2.

\item[\emph{(iii)}] $\operatorname{div}\alpha=nA-nB$. Then $A$ is a smooth point of $\mathcal X$ if and only if $P(0)\ne 0$

\item[\emph{(iv)}] Let $(x_0,\alpha_0)$ be a point on $\mathcal X$, $\alpha_0\ne 0$. It is a critical point of $\alpha$ if $x_0$ satisfies the equation $$n^2 x_0^{n-1}-nP(x_0)P'(x_0)+x_0(P'(x_0))^2=0$$
\end{itemize}
\end{utv}

\textbf{Proof:}

\emph{(i)} Let's specify the coordinates $(x, y)$ of $\mathcal X$ so that its equation has the form $y^2=f(x)$ and $x(B)=\infty$, so the Abel function takes the form $\alpha= P (x) + Q (x) y $. Then the $\mathcal X$'s equation takes the form $ \alpha ^ 2 + 2\alpha P (x) + x^n = 0 $. $ \Box$

\emph{(ii)} $\mathcal X$ is elliptic, $\alpha= P (x) + Q (x) y $ and $\deg \alpha=n$. So $\deg P\le\frac n2$. For even $n$ we have: $P^2/4-x^n=f(x)Q^2(x)$, $\deg f=3$, so degree of $P^2/4-x^n$ is odd, and the leading coefficient of $P^2$ is $4$. $\Box$

\emph{(iii)} $\frac{\partial F}{\partial x}=\alpha+nx^{n-1}$ and $\frac{\partial F}{\partial \alpha}=2\alpha+P(x)$. $x(A)=0$ and $\alpha=0$, so $A$ is singular if and only if $P(0)=0$. $ \Box $

\emph{(iv)} $\alpha_0$ is a critical value of $\alpha$ if and only if $F(x,\alpha_0)$ has a multiple root, i.e.  $F(x_0,\alpha_0)=0$ and $\left.\frac{\partial  F(x,\alpha)}{\partial x}\right|_{(x_0,\alpha_0)}=0$. But $\left.\frac{\partial F(x,\alpha)}{\partial x}\right|_{(x_0,\alpha_0)}=\alpha_0 P'(x_0)+nx_0^{n-1}=0$, i.e. $\alpha_0=-\frac{nx_0^{n-1}}{P'(x_0)}$, so $x_0$ is a root of $F(x,-\frac{n x^{n-1}}{P'(x)})$. $\Box$

\begin{notation}
We introduce the polynomials $$R_n(x):=-\left(\frac{-1+\sqrt{1-4x}}2\right)^n-\left(\frac{-1-\sqrt{1-4x}}2\right)^n=-2^{1-\frac{n}2}x^{\frac{n}2}T_{n}\left(-\frac1{2\sqrt{x}}\right)$$ where $T_n$ is a Chebyshev polynomials. For example $R_0(x)=-2$, $R_1(x)=1$, $R_2(x)=2x-1$, $R_3(x)=-3x+1$, $R_4(x)=-2x^2+4x-1$, $R_5(x)=1-5x+5x^2,\ldots$
\end{notation}

\begin{utv}
\begin{itemize}
\item[\emph{(i)}]  Let $\mathcal X$ be an algebraic curve over $\Bbbk$ defined by the equation $\alpha^2+\alpha P(x)+x^n=0$, $\deg P\le n/2$ and suppose $\operatorname{char}\Bbbk\nmid n$. The curve $\mathcal X$ is rational if and only if $P(x)=C^{2k-n} x^{k}R_{n-2k}(C^2 x)$, where $k\le n/2$ and $C\in\Bbbk$.

\item[\emph{(ii)}] Let the pair $(\mathcal X,\alpha)$ be defined by the equation $\alpha^2+\alpha x^k R_{n-2k}(x)+x^n=0$. Then there exists $t\in\Bbbk(\mathcal X)$ such that $\alpha=t^{n-k}(-t-1)^k$.
\end{itemize}
\end{utv}

\textbf{Proof:} (i) First, let $P(0)\ne 0$, so $A$ is a smooth point of $\mathcal X$. In the notation of the previous proof we obtain $\frac{P^2(x)}4-x^n=f(x)Q^2(x)$. $\mathcal X$ is rational if and only if $\deg f\le 2$. Besides, from $\deg (\frac{P^2}4-x^n)=n$ we get $\deg f=1$ for odd $n$ and $\deg f=2$ for even $n$.

Next we make the change of variables in the identity $\frac{P^2}4-x^n=f(x)Q^2(x)$; let $x:=\frac 1{t^2}$.
So $(t^n\cdot P(\frac 1{t^2})/2)^2=f(\frac 1{t^2})Q^2(\frac 1{t^2})\cdot t^{2n}+1$.
From this identity we see that $t^n\cdot P(\frac 1{t^2})$ is a Shabat polynomial and get valencies of its critical point. So we get $t^n\cdot P(\frac 1{t^2})=2T_n(Const\cdot t)$.

Secondly, let us consider the case $P(0)=0$. We have $P(x)=x^kP_1(x)$, $P_1(0)\ne 0$. We reduce this case to the previous one by the change of variables $\alpha_1=\alpha/x^k$.

(ii) $t=\frac{-1+\sqrt{1-4x}}2=-\frac{2\alpha+2Q+P}{4Q}\in\Bbbk(\mathcal X)$. Then $t(-t-1)=x$ and $t^{m}+(-1-t)^m=R_m(x)$. So $\alpha=x^k t^{n-2k}=t^{n-k}(-t-1)^k$. $\Box$

\section{Abel pairs of genus 1}

\subsection{The family of Abel pairs of genus 1}

The modular curve $Y_1(n)$ (see ~\cite{Firstcourse}) parameterizes the pairs $(\mathcal E,A-B)$, where $\mathcal E$ is an elliptic curve and $A-B$ is a divisor on it of order exactly $n$. The modular curve $X_1(n)$ compactifies $Y_1(n)$, i.e. over the punctures of $Y_1(n)$ the stable curves are added.

An Abel pair of genus 1 is actually determined by arbitrary elliptic curve and two points $A$, $B$, such that $n(A-B)\equiv 0$. But the true order of divisor $(A-B)$ may be less than $n$; it means that such an Abel function is a power of another Abel function.

Next we will see $X_1(n)$ is the base of family of Abel pairs.

\begin{definition}
An Abel pair $(\mathcal X,\alpha)$ is called \emph{imprimitive} if there exists an other Abel pair $(\mathcal X,\alpha_0)$ and natural $k>1$ such that $\alpha=\alpha_0^k$. Otherwise, an Abel pair is called \emph{primitive}.
\end{definition}

\begin{theorem}\label{components}
\begin{itemize}
\item[(i)] The parameter space of Abel pairs of genus 1 and degree $n$ has $\sigma_0(n)-1$ components, where $\sigma_0(n)$ is the number of divisors of n.

\item[(ii)] Number $d$ such that $d\mid n$ and $1<d$ corresponds to component of parameter space of Abel pairs of genus 1 and degree $n$; it consists of Abel pairs which are the $\frac n d$-powers of primitive Abel pair of genus 1.

\item[(iii)] $Y(n)$ is isomorphic to the space of parameters of elliptic primitive Abel pairs.
\end{itemize}
\end{theorem}

\textbf{Proof:} See ~\cite{Chen}, proposition 3.2. $\Box$

\begin{theorem}\label{kappabelaja}
\begin{itemize}
\item[\emph{(i)}] Besides $0$ and $\infty$, a generic Abel function $\alpha$ (except the finite set of Abel functions) on a curve of genus 1 has exactly two more critical values: $\{k_1,k_2\}$.

\item[\emph{(ii)}] $\varkappa_n=\cfrac{2-\cfrac{k_1}{k_2}-\cfrac{k_2}{k_1}}4=-\cfrac{(k_1-k_2)^2}{4k_1k_2}$ is a well-defined Belyi function on $X_1(n)$.
\end{itemize}
\end{theorem}

\textbf{Proof:}

\emph{(i)} Let $e_P$ be a ramification index of $\alpha$ at the point $P$. Then by the Riemann-Hurwitz formula we get $2\deg \alpha=\sum\limits_{P\in\mathcal X}(e_P-1)$. But $\operatorname{div}(\alpha)=nA-nB$, hence $e_A=e_B=n=\deg \alpha$.

Consequently $2=\sum\limits_{P\in\mathcal X\backslash\{A,B\}}(e_P-1)$. Therefore, we have two possibilities: two points with ramification index 2 (denote them $C_1$ and $C_2$) or one point with ramification index 3. So $\alpha$ has one or two critical values besides $0$ and $\infty$. The case of one additional critical value corresponds to the case of Abel-Belyi pairs in which we assume that $k_1=k_2$. Abel-Belyi pairs with degree $n$ constitute a finite set.

\emph{(ii)} $\alpha(C_1)=k_1$ and $\alpha(C_2)=k_2$ are defined up to multiplication by a constant, common inversion and permutation. So $\cfrac{k_1}{k_2}+\cfrac{k_2}{k_1}$ is a well-defined function on $X_1(n)$, and $\varkappa_n$ is a well-defined function too.

Denote $\lambda(t)=\frac{2-t-1/t}4$, so $\varkappa_n=\lambda(k_1/k_2)$. The critical values of $\cfrac{2-\cfrac{k_1}{k_2}-\cfrac{k_2}{k_1}}4$ are either the critical points of $\lambda(t)$ or the values of $\lambda$ in the critical values of $k_1/k_2$. The critical values of $\lambda$ are $0$, $1$, $\infty$. The critical values of $k_1/k_2$ are $0$, $1$ and $\infty$, because $k_1/k_2$ is a local parameter in any other point. So $\varkappa_n$ is a Belyi function. $\Box$

\subsection{Abel-Belyi pairs of genus 1}

\begin{utv}
Dessins on the torus corresponding to the Abel-Belyi pairs have sets of valencies
$ (n | n | 3,1,1, \ldots, 1) $ or $ (n | n | 2,2,1, \ldots, 1) $.
\end{utv}

\textbf{Proof:} For a toric dessin  corresponding with the Abel-Belyi pair with the set of valencies $ (n | n | a_1, \ldots, a_k) $, we have by Riemann–Hurwitz formula $ k = n-2 $. Since $ a_1 + \ldots + a_k = n $, we obtain the desired assertion. $ \Box$

\begin{utv}
\emph{(i)} A dessin with a set of valencies $ (n | n | 3,1, \ldots, 1) $ is determined uniquely by the set $ \langle a, b, c \rangle $, defined up to a cyclic permutation, where $ a + b + c = n $.

We denote this dessin by $\hexagon_{a,b,c}$
\begin{figure}[h]
\begin{center}
\includegraphics[width=4in,height=1.3in]{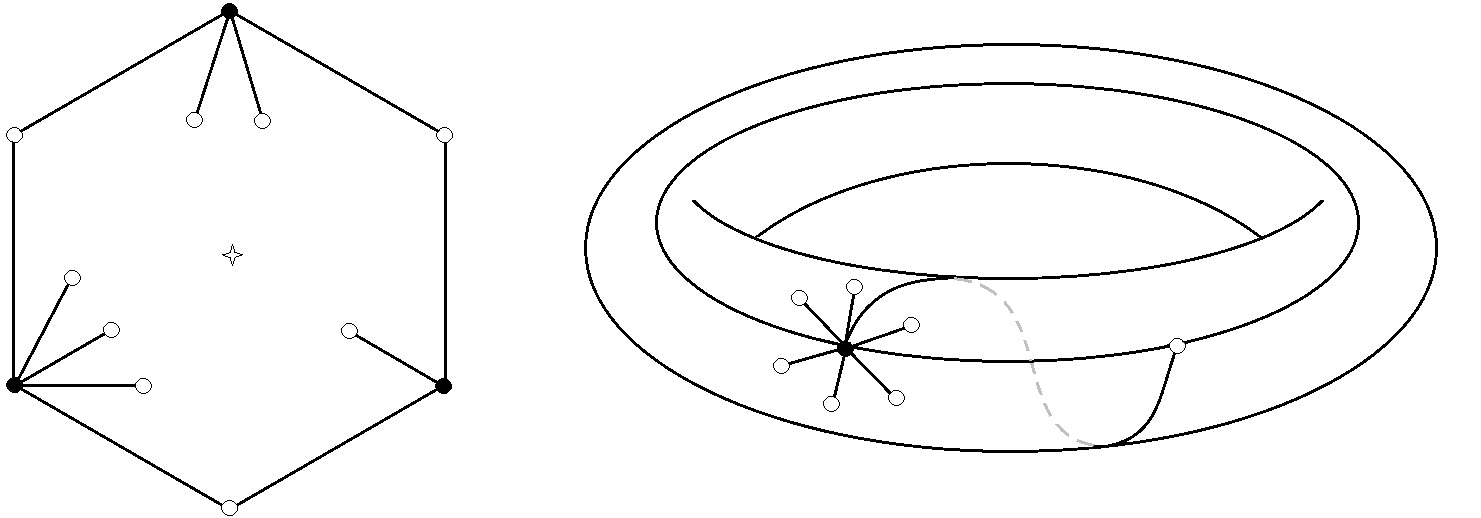}
\caption{$\hexagon_{4,2,3}$}
\end{center}
\end{figure}

\emph{(ii)} A dessin with a set of valencies $ (n | n | 2,2, \ldots, 1) $ is uniquely defined by a set $\langle a, b, c, d \rangle $, defined up to a cyclic permutation.

We denote this dessin by $\Box_{a,b,c,d}$

\begin{figure}[h]
\begin{center}
\includegraphics[width=4in,height=1.3in]{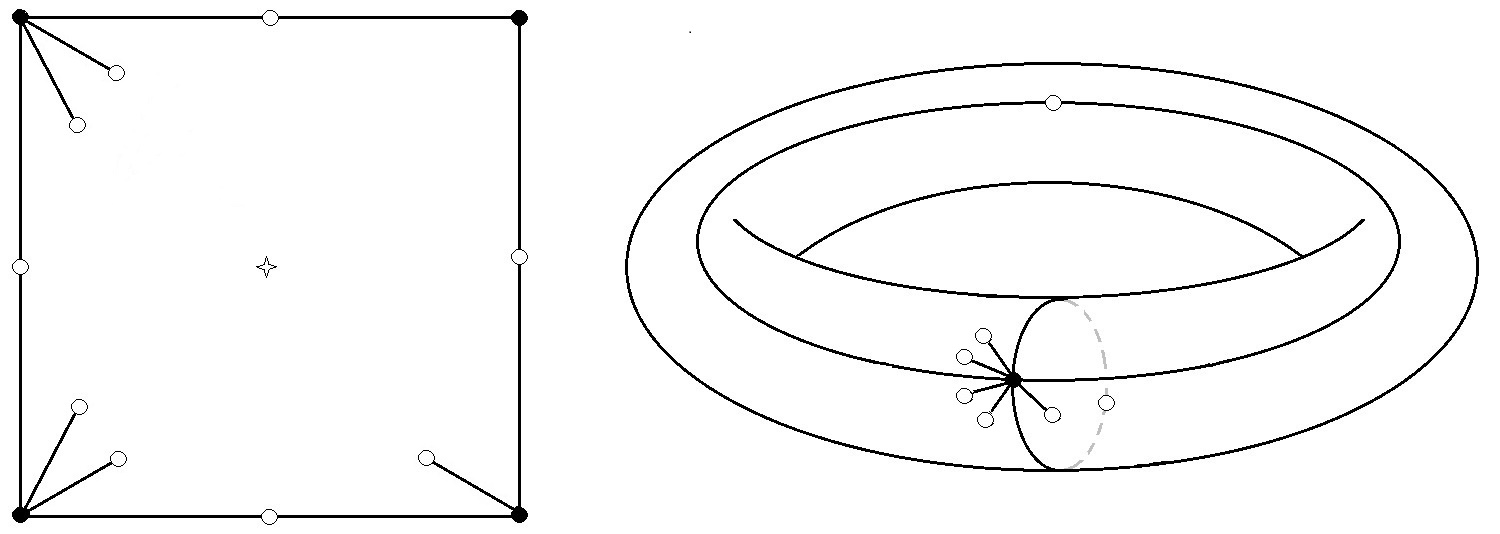}
\caption{$\Box_{3,2,1,3}$}
\end{center}
\end{figure}

\end{utv}

\textbf{Proof:} (i) Let's consider a dessin with a set of valencies $ (n | n | 2,2, \ldots, 1) $. This toric dessin has one white vertex of valence 3, the others have valencies 1, i.e. they are terminal vertices. Let's erase the terminal vertexes and edges, than the valency list takes the form $(3|3|3)$.

Next we return to the dessin $(3|3|3)$ previously erased terminal white vertices. They can be added at one of three angles between 3 edges emanating from the black vertex of $(3|3|3)$. Let the number of these edges read contrary-clockwise be $ a-1 $, $ b-1 $, $ c-1 $. This numbers define the original dessin $\hexagon_{a,b,c}$.

The case of $\Box_{a,b,c,d}$ is treated in a similar way. $\Box$

\begin{theorem}\label{degtheorem}
\begin{itemize}
\item[\emph{(i)}]
A toric dessin $ \hexagon_{a,b,c}$ corresponds to a Belyi function, which is an
$ m $-th power if and only if $ m $ divides $ \operatorname{gcd}(a, b, c) $.

\item[\emph {(ii)}] A toric dessin $\Box_{a,b,c,d}$ corresponds to a Belyi function,
which is an $ m $-th power if and only if $ a \equiv -b \equiv c \equiv -d \; \; (mod \;m) $.
\end{itemize}
\end{theorem}

\textbf{proof:} See the theorem 3 in \cite{Pakovich}. $\Box$

\begin{notation}
Let's denote $m_{\hexagon}(n,\Bbbk)$ and $m_{\Box}(n,\Bbbk)$ the numbers of Abel-Belyi pairs over $\Bbbk$ with set of valencies $(n | n | 3,1, \ldots, 1)$ and $(n | n | 2,2,1 \ldots, 1)$ respectively. Let $\widetilde{m}_{\hexagon}(n,\Bbbk)$ and $\widetilde{m}_{\Box}(n,\Bbbk)$ denote the numbers of such primitive pairs.

Recall that we count Abel-Belyi pairs $(\mathcal X,\alpha)$ with the weight $\frac 1{|\operatorname{Aut}(\mathcal X,\alpha)|}$.

So $m_{\hexagon}(6,\mathbb C)=3\frac13$. Indeed here is the list of Abel-Belyi pairs over $\mathbb C$ with the set of valencies $(6|6 |3,1,1,1)$: $\hexagon_{4,1,1}$, $\hexagon_{3,2,1}$, $\hexagon_{3,1,2}$ and $\hexagon_{2,2,2}$, the last one with the automorphism group of order 3. But $\widetilde{m}_{\hexagon}(6,\mathbb C)=3$.

Next, $m_{\Box}(6,\mathbb C)=2\frac12$. The list of Abel-Belyi pairs over $\mathbb C$ with set of valencies $(6|6 |2,2,1,1)$ is $\Box_{2,2,1,1}$, $\Box_{3,1,1,1}$, $\Box_{2,1,2,1}$, the last one with the automorphism group of order 2. But $\widetilde{m}_{\Box}(6,\mathbb C)=1$, and only $\Box_{2,2,1,1}$ is primitive. $\Box_{3,1,1,1}$ corresponds to the Belyi function, which is a square and $\Box_{2,1,2,1}$ corresponds to a cube.  \end{notation}

\begin{theorem}\label{combinatoric}
If $n>3$ then
\begin{itemize}
\item[\emph{(i)}] $m_{\hexagon}(n,\mathbb C)=\frac{(n-1)(n-2)}6;\,\,\,\,\,\,\, m_{\Box}(n,\mathbb C)=\frac{(n-1)(n-2)(n-3)}{24};$
\item[\emph{(ii)}] $m_{\hexagon}(n,\mathbb C)=\sum\limits_{1<d, d\mid n}\widetilde{m}_{\hexagon}(d,\mathbb C);$
\item[\emph{(iii)}] $\widetilde{m}_{\hexagon}(n,\mathbb C)= \cfrac {\varphi (n) \psi (n)}6 - \cfrac {\varphi (n)}2;\,\,\,\,\,\,\,\widetilde{m}_{\Box}(n,\mathbb C)=\cfrac{(n-6)\varphi(n)\psi(n)}{24}+\cfrac{\varphi(n)}2;$
where $\varphi$ is the Euler function and $\psi $ is the Dedekind psi function $\psi(n) = n \prod \limits_ {p | n} \left( 1 + \cfrac 1p \right) $.
\end{itemize}
\end{theorem}

\textbf{Proof:} (i) The number of $\hexagon_{a,b,c}$ is the number of solutions of the equation $a+b+c=n$ in positive integers divided by 3 (because of rotations of the hexagon). So $m_{\hexagon}(n,\mathbb C)={n-1\choose 3}/3$. Similarly, $m_{\Box}(n,\mathbb C)={n-1\choose 4}/4$.

(ii) According to ~\ref{degtheorem}, an Abel-Belyi pair $\hexagon_{a,b,c}$ is primitive if and only if $\operatorname{gcd}(a,b,c)=1$. Also, for all $a,b,c,k$, $\alpha_{\hexagon{ka,kb,kc}}=\alpha_{\hexagon{a,b,c}}^k$. All the Abel-Belyi functions $\alpha_{\hexagon{a,b,c}}$ are powers of the primitive Abel-Belyi functions. So the number $m_{\hexagon}(n,\mathbb C)$ of all pairs is the sum of the primitive ones with degree $d\mid n$, i.e. $\widetilde{m}_{\hexagon}(d,\mathbb C)$.

(iii) Apply to (ii) the M\"{o}bius inversion formula $\widetilde{m}_{\hexagon}(n,\mathbb C)=\sum\limits_{d\mid n}m_{\hexagon}(d)\mu(n/d)$. By using (i) we get $\widetilde{m}_{\hexagon}(n,\mathbb C)=\sum\limits_{d\mid n} \frac{(d-1)(d-2)}6\mu(n/d)= \frac {\varphi (n) \psi (n)}6 - \frac {\varphi (n)}2$.

The formula for $m_{\Box}(n,\mathbb C)$ is proved in the same way by more cumbersome considerations. $\Box$

\subsection{Dessin on $X_1(n)$ corresponding to $\varkappa_n$}

In the next proposition we classify the critical points of $\varkappa_n$.

\begin{theorem}\label{kappazeros}
Let $\operatorname{char}\Bbbk \neq 2,3$, $\operatorname{char}\Bbbk\nmid n$. The critical points of $\varkappa_n$ on $Y_1(n)$ belong to one of these three types:
\begin{itemize}
\item[\emph{(i)}] Each primitive Abel-Belyi pair of the type $\hexagon_{a,b,c}$ corresponds to a zero of multiplicity 3 of $\varkappa_n$;

\item[\emph{(ii)}] Each primitive Abel-Belyi pair of the type $\Box_{a,b,c}$ corresponds to a zero of multiplicity 2 of $\varkappa_n$;

\item[\emph{(iii)}] Other critical points of $\varkappa_n$ correspond to the critical value $\varkappa_n=1$ and have multiplicity 2.
\end{itemize}
\end{theorem}

\textbf{Proof:} $\varkappa_n=\cfrac{2-\cfrac{k_1}{k_2}-\cfrac{k_2}{k_1}}4=-\frac{(k_1-k_2)^2}{4k_1k_2}$. The function $\frac1{4k_1k_2}$ has no zeros or poles on $Y_1(n)$ so the zeros are defined only by $(k_1-k_2)^2$ zeros. This function has its zeros on $Y_1(n)$ if $k_1=k_2$, i.e in the case of Abel-Belyi pairs.

Now we consider primitive pairs $\hexagon_{a,b,c}$ and $\Box_{a,b,c}$. In the notations of ~\ref{kappabelaja} $(k_1-k_2)^2=(\alpha(C_1)-\alpha(C_2))^2$.

In the case $\Box_{a,b,c,d}$ critical points $C_1$ and $C_2$ are different. Locally $\alpha(C_1)$ and $\alpha (C_2)$ are well-defined functions on the base and $\alpha(C_1)-\alpha(C_2)$ has a simple zero because it is a local parameter ($\operatorname{char}\Bbbk \neq 2,3$). So $\varkappa_n$ has a zero of multiplicity 2 in $\Box_{a,b,c,d}$.

In the case $\hexagon_{a,b,c}$ the critical points $C_1$ and $C_2$ coincide. Then $\frac{\alpha(C_1)-\alpha(C_2)}{x(C_1)-x(C_2)}$ and $x(C_1)-x(C_2)$ are local parameter on the base and have a simple zero in $\hexagon_{a,b,c}\in X_1(n)$. So in $\varkappa_n$ has a zero of multiplicity 3 in $\hexagon_{a,b,c}$.

By ~\ref{kappabelaja} the other critical points of $\varkappa_n$ correspond to the case $\varkappa_n=1$. In that point $k_1=-k_2$ and $k_1/k_2$ is a local parameter, so $\varkappa_n=\frac{2-\frac{k_1}{k_2}-\frac{k_2}{k_1}}4=1-\frac{(k_1+k_2)^2}{4k_1k_2}$ has multiplicity 2. $\Box$

\begin{sled}\label{kappadiv}
Let $\hexagon^n_1$, $\hexagon^n_2,\ldots$ be the list of primite Abel-Belyi pairs with sets of valencies $(n | n | 3,1,1, \ldots, 1)$, and $\Box^n_1$, $\Box^n_2,\ldots$ be the list of primitive Abel-Belyi pairs with sets of valencies $(n | n | 2,2,1, \ldots, 1)$, then $$\operatorname{div}(\varkappa_n)=3\cdot(\hexagon^n_1+\hexagon^n_2+\ldots)+2(\Box^n_1+\Box^n_2+\ldots)+(points \,from\, X_1(n)\backslash Y_1(n))$$
\end{sled}

\begin{figure}[h]
\begin{center}
\includegraphics[width=3.8in,height=1.7in]{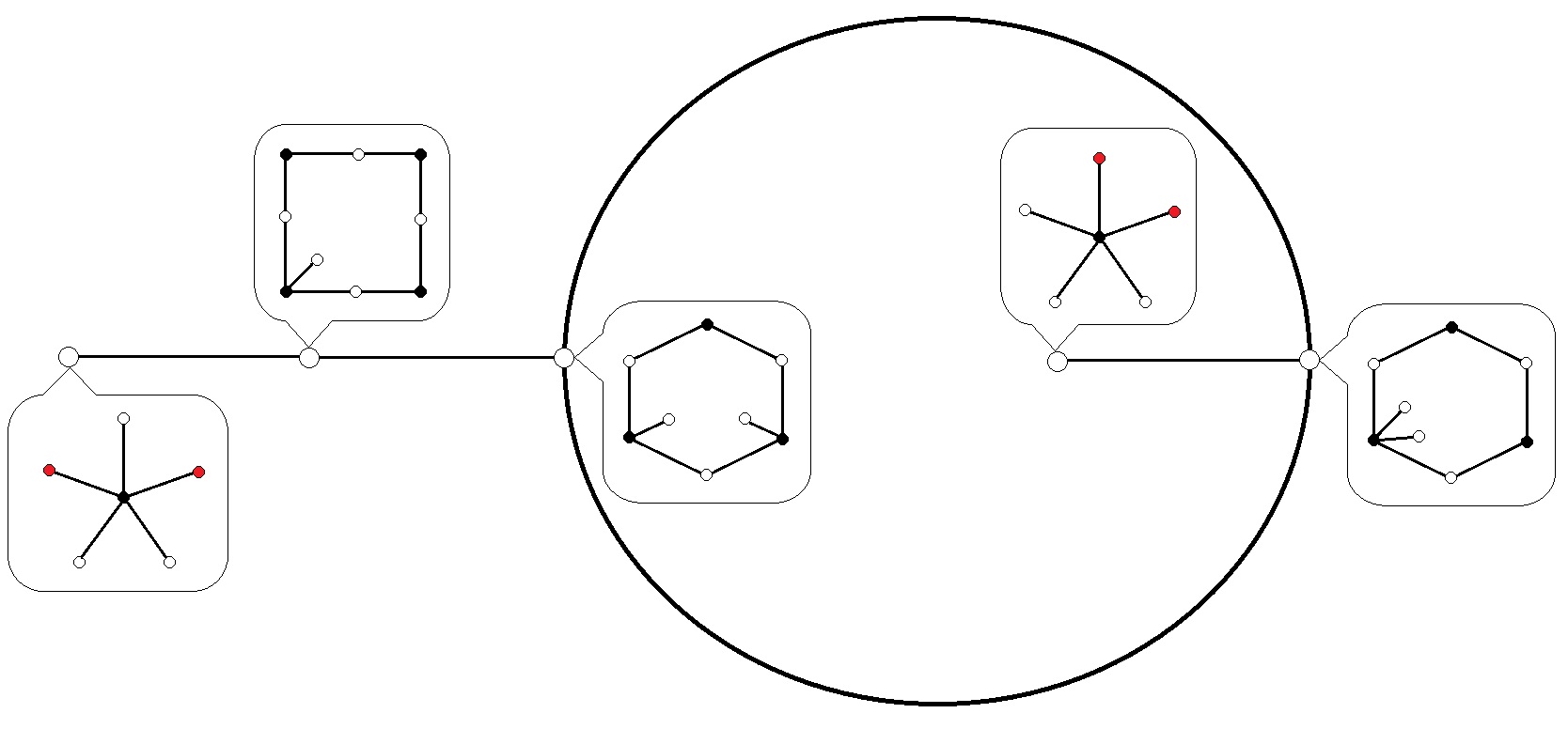}
\caption{Dessin on $X_1(5)$ corresponding to $\varkappa_5$}
\end{center}
\end{figure}

\begin{theorem}\label{kappapoles}
Let $\operatorname{char}\Bbbk \neq 2,3$, $\operatorname{char}\Bbbk\nmid n$. Any point of $X_1(n)\backslash Y_1(n)$ belongs to one of the following types.

\begin{itemize}
\item[\emph{(i)}] Points whose Abel-Belyi pairs have representatives $\alpha^2+\alpha R_n(x)+x^n=0$. $X_1(n)\backslash Y_1(n)$ contains $\frac{\varphi(n)}2$ such points. In these points $\varkappa_n$ have simple zeros.

\item[\emph{(ii)}] Points whose Abel-Belyi pairs have representatives $\alpha^2+\alpha x^k R_{n-2k}(x)+x^n=0$, where $k$ satisfy next conditions: $0<k<\frac n2$, $\operatorname{char}\Bbbk\nmid k$, $\operatorname{char}\Bbbk\nmid (n-k)$. $X_1(n)\backslash Y_1(n)$ contains $\varphi(\operatorname{gcd}(n,k))$ such points. In these points $\varkappa_n$ have poles of multiplicity $\frac{k(n-k)}{\operatorname{gcd}(n,k)}$.

\item[\emph{(iii)}] Points whose Abel-Belyi pairs have representatives $\alpha^2+2\alpha x^{n/2}+x^n=0$. This case is possible only for even $n$. $X_1(n)\backslash Y_1(n)$ contains $\frac{\varphi(n/2)}2$ such points. In these points $\varkappa_n$ have poles of multiplicity $n/2$.

\item[\emph{(iv)}] Points whose Abel-Belyi pairs have representatives  $\alpha^2+\alpha x^k R_{n-2k}(x)+x^n=0$, where $\operatorname{char}\Bbbk>0$ and $k$ satisfy the conditions $0<k< n$, $\operatorname{char}\Bbbk\mid k$. $X_1(n)\backslash Y_1(n)$ contains $\varphi(\operatorname{gcd}(n,k))$ such points. In these points $\varkappa_n=1$ if $k$ is odd and $\varkappa_n=0$ if $k$ is even. The multiplicity of $\varkappa_n$ in these critical points is $\frac{k (n-k)}{p^{\operatorname{ord}_p(k)}\cdot\operatorname{gcd}(n,k)}$.
\end{itemize}
\end{theorem}

\textbf{Proof:} In the case (i) the Abel pair in the notations of \ref{kappabelaja} correspond to $k_1=k_2=1$ and $C_1=C_2$, but this point is singular (double point). By 4.4(ii) after the resolution of singularity this Abel pair takes the form $(\mathbb P^1(\Bbbk),t^n)$; in this form the double point has two preimages, in which $t$ assumes the values belonging to $\sqrt[n]{1}$. To count such pairs we should choose two different points with $t\in\sqrt[n]{1}$ as the normalization of preimages of $C$. Due to the rotations we have $n/2$ Abel pairs. Excluding Abel pair that lie over $X_1(d)$, $d\mid n$, $d<n$, we have $\varphi(n)/2$ pairs.

(ii)-(iii) In this case we have a pole of $\varkappa_n$ on $X_1(n)\backslash Y_1(n)$, i.e. $k_1k_2=0$.

We use 4.6(ii), $\alpha=t^{n-k}(t-1)^k$.
To calculate multiplicity of this pole let's consider the locally critical value $k_1$ going to zero. We denote $x_1=x(C_1)$. Locally at $t=0$ and $t=1$, $\alpha$ has the form $At^k$, $B(t-1)^{n-k}$ for some $A$, $B\in\Bbbk$. By changing variables back to $x$ we get
$k_1=Ax_1^{1/k}=Bx_1^{1/(n-k)}$. This equation has $\frac{k(n-k)}{\operatorname{gcd}(n,k)}$ solutions by $x_1$ because $\operatorname{char}\Bbbk\nmid k$, $\operatorname{char}\Bbbk\nmid (n-k)$. So the multiplicity of this pole of $\varkappa_n$ is $\frac{k(n-k)}{\operatorname{gcd}(n,k)}$ too. Calculation of the number of such poles is the same as in (i)

(iv) In this case both $k_1$ and $k_2$ go to $0$.
We want to calculate the undeterminate form $k_1/k_2$. We denote $x_1=x(C_1)$, $x_2=x(C_2)$. Locally we have $x_1=-x_2$. So $k_2/k_1=(-1)^k$.
I.e. in thise points $\varkappa_n=1$ if $k$ is odd and $\varkappa_n=0$ if $k$ is even.
Multiplicity and calculation of the number of such points is same as in (i), (ii) $\Box$

\begin{utv}\label{kappadeg}
$\deg{\varkappa_n}=\frac{ n\varphi(n)\phi(n)}{12}$.
\end{utv}

\textbf{Proof:} From the previous theorem  $\deg{\varkappa_n}=\frac12\sum\limits_{k=1}^{n-1}\frac{\varphi(\operatorname{gcd}(n,k))}{\operatorname{gcd}(n,k)}k(n-k)=\frac{ n\varphi(n)\psi(n)}{12}$. $\Box$

To calculate the numbers $m_{\Box}$ and $m_{\hexagon}$ we need also to take into account the imprimitive Abel pairs. So we consider all the components of bases of families of the Abel pairs, i.e. $X_1(d)$'s for $1<d$, $d\mid n$. $\frac{2-\frac{k_1}{k_2}-\frac{k_2}{k_1}}4$ is a well-defined Belyi function on all these components . We denote it $\varkappa_{d, n}$.

\begin{utv}\label{kappapower}
 If $d>1$ is a divisor of $n$, then $\varkappa_{d, n}=\frac{T_{n/d}\circ \varkappa_d+1}2$.
\end{utv}

\textbf{Proof:} By ~\ref{components} $\varkappa_{d, n}=\cfrac{2-\cfrac{k_1}{k_2}-\cfrac{k_2}{k_1}}4=\cfrac{2-\cfrac{\kappa_1^{n/d}}{\kappa_2^{n/d}}-\cfrac{\kappa_2^{n/d}}{\kappa_1^{n/d}}}4=\frac{T_{n/d}(\varkappa_d)+1}2$,
where $\kappa_1$ and $\kappa_2$ are critical points of an imprimitive Abel pair of genus 1 and degree $n/d$. $\Box$

\section{Number of Abel-Belyi pairs of genus 1 }

Excluding statements concerning the lists of Abel-Belyi pairs of genus 1, obtained by the topological claims, all the statements above are valid in the case of positive characteristic (under the condition $\operatorname{char}\Bbbk \nmid n$, $\operatorname{char}\Bbbk\ne 2,3$).

Next we count the number of Abel-Belyi pairs of genus 1 in positive characteristic. Instead of topological considerations which not available now, we use the algebraic ones: especially, study of the divisor of $\varkappa_n$. We will use ~\ref{kappazeros} -- ~\ref{kappapower}.

\begin{utv}
\begin{itemize}
\item[\emph{(i)}] $2m_{\Box}(n,\mathbb C)+3m_{\hexagon}(n,\mathbb C)=\frac{(n-1)(n-2)(n+3)}{12}$

\item[\emph{(ii)}] $\deg\varkappa_n(\mathbb C)-\deg\varkappa_n(\overline{\mathbb F_p})=\sum\limits_{k=1}^{\lfloor n/p\rfloor}\frac{\varphi(\operatorname{gcd}(k,n))}{\operatorname{gcd}(k,n)}kp(n-kp).$

\item[\emph{(iii)}] $2m_{\Box}(n,\mathbb C)+3m_{\hexagon}(n,\mathbb C)-2m_{\Box}(n,\overline{\mathbb F_p})-3m_{\hexagon}(n,\overline{\mathbb F_p})=\sum\limits_{k=1}^{\lfloor n/p\rfloor}kp(n-kp)+\sum\limits_{k=1}^{\lfloor n/2p\rfloor}{\frac {2k}{p^{\operatorname{ord}_p(k)}}(n-2kp)}$
\end{itemize}
\end{utv}

\textbf{Proof:} For (i) use ~\ref{kappapower}: $\sum\limits_{1<d\mid n}\deg{\varkappa_{d, n}}=\sum\limits_{1<d\mid n}\frac nd\deg{\varkappa_{d}}=n\sum\limits_{1<d\mid n}\frac{ \varphi(d)\psi(d)}{12}=\frac{n(n-1)(n+1)}{12}$.

By ~\ref{kappazeros} -- ~\ref{kappapoles} the function $\varkappa_{d, n}$ has $\varphi(d)/2$ simple zeros described in ~\ref{kappapoles}(i) besides the points of $X_1(d)$ over which Abel-Belyi pairs of genus 1 lie. By summing the numbers of all the simple zeros of $\varkappa_{d, n}$, we have $\sum\limits_{1<d,\, d\mid n}\frac{\varphi(d)}2=\frac{n-1}2$ simple zeros. So from ~\ref{kappadiv} we get $2m_{\Box}(n)+3m_{\hexagon}(n)=\frac{n(n-1)(n+1)}{12}-\frac{n-1}2=\frac{(n-1)(n-2)(n+3)}{12}$.

For (ii) we calculate the difference between the degrees of $\varkappa_n(\mathbb C)$ and $\varkappa_n(\overline{\mathbb F_p})$. Let us compare the list of poles of $\varkappa_n$ over $\mathbb C$ and $\overline{\mathbb F_p}$. By ~\ref{kappapoles} (iii) in the case of $\overline{\mathbb F_p}$ the poles of $\varkappa_n$ of degree $\frac{kp(n-kp)}{\operatorname{gcd}(k,n)}$ with $0<k<\frac np$ disappear.

For (iii) by ~\ref{kappapoles} (iv) in characteristic $\operatorname{char}\Bbbk=p$ the list of zeros of $\varkappa_n(\overline{\mathbb{F}_p})$ as compared with the case of $\operatorname{char}\Bbbk=0$ has additional zeros on $X_1(n)\backslash Y_1(n)$. List of its multiplicities is $\frac{2k(n-2kp)}{p^{\operatorname{ord}_p(k)}\operatorname{gcd}(n,k)}$, $0<k<\frac n{2p}$. By summing this multiplicities of zeros for $\varkappa_{d, n}(\overline{\mathbb{F}_p})$, ($1<d$, $d\mid n$) and by (ii) we calculate the difference between the number of zeros of $\varkappa_n(\mathbb C)$ and $\varkappa_n(\overline{\mathbb{F}_p})$ on $Y_1(n)$. By ~\ref{kappazeros} this is $2m_{\Box}(n,\mathbb C)+3m_{\hexagon}(n,\mathbb C)-2m_{\Box}(n,\overline{\mathbb F_p})-3m_{\hexagon}(n,\overline{\mathbb F_p})$. $\Box$

The previous theorem is insufficient to calculate $m_{\Box}(n,\overline{\mathbb F_p})$ and $m_{\hexagon}(n,\overline{\mathbb F_p})$.
To calculate these quantities separately, we need one more equation. We will use the genus of $X_1(n)$.

\begin{lemma}
$genus(X_1(n,\overline{\mathbb F_p}))=genus(X_1(n))=\frac{\varphi(n)\psi(n)}{24}-\frac{\varphi(n)\sigma_0(n)}{4}+1$;
\end{lemma}

\textbf{Proof:} For $\operatorname{char}\Bbbk=0$ we can apply ~\ref{combinatoric}, ~\ref{kappazeros}, ~\ref{kappapoles}(i)--(iii) and the Riemann–Hurwitz formula  to the function $\varkappa_n$ on $X_1(n)$. By Igusa's theorem for $p\nmid n$ the curve $X_1(n)$ has a good reduction at $p$. See  8.6.1. at ~\cite{Firstcourse}. $\Box$

\begin{theorem}
\begin{itemize}
\item[\emph{(i)}] $m_{\hexagon}(n,\mathbb C)-m_{\hexagon}(n,\overline{\mathbb F_p})=\sum\limits_{0<k<n/p} (n-kp)$

\item[\emph{(ii)}] $m_{\Box}(n,\mathbb C)-m_{\Box}(n,\overline{\mathbb F_p})=\sum\limits_{k=1}^{\lfloor n/p\rfloor}\frac{kp-3}2(n-kp)+\sum\limits_{k=1}^{\lfloor n/2p\rfloor}{\frac {k}{p^{\operatorname{ord}_p(k)}}(n-2kp)}$
\end{itemize}
\end{theorem}

\textbf{Proof:} We apply Riemann–Hurwitz formula to $\varkappa_n(\mathbb C)$ and $\varkappa_n(\overline{\mathbb F_p})$ to calculate $genus(X_1(n,\overline{\mathbb F_p}))$ and $genus(X_1(n))$. From the previous lemma $genus(X_1(n,\overline{\mathbb F_p}))=genus(X_1(n,\mathbb C))$, so we consider the difference of these two formulas. We also sum up these differences over all the divisors of $n$ (similarly to 6.11), and we get $2(m_{\hexagon}(n)-m_{\hexagon}(n,\overline{\mathbb F_p}))+(m_{\Box}(n)-m_{\Box}(n,\overline{\mathbb F_p}))=\frac{p+1}2(n-p)+(p+1)(n-2p)+\frac{3p+1}2(n-3p)+2(p+1)(n-4p)+\ldots$

With 6.11 we get (ii) and $m_{\hexagon}(n,\overline{\mathbb F_p})=m_{\hexagon}(n, \mathbb C)-(n-p)-(n-2p)-\ldots$. $\Box$


\begin{zamech}
Note that in the case $p>n$ we have $m_{\Box}(n,\overline{\mathbb F_p})=m_{\Box}(n,\mathbb C)$ and $m_{\hexagon}(n,\overline{\mathbb F_p})=m_{\hexagon}(n,\mathbb C)$.
\end{zamech}

\section{Example}

\subsection{Method: Pad\'e approximation}

We describe our method for calculating the Abel pair of genus 1:

\begin{definition}
\emph{Pad\'e approximation} (see ~\cite{pade}) of order $[n, m]$ of a real-valued function $f(x)\in C^{n+m}(U)$, (where $U\subset\mathbb R$ is a neighborhood of $0$) is the ratio of two polynomials $R_{[n, m]} = \frac {p_{[n, m]} (x)}{q_{[n, m]} (x)} $, $ \deg p_ {[n, m] } (x) \le m $, $ \deg q_ {[n, m]} (x) \le n $,
for which $ f^{(i)} (0) = R^{(i)}_{[n, m]} (0) $ for $ 0 \le i \le m + n $.
\end{definition}

Let $\mathcal E$ be an elliptic curve. We want to find $\alpha\in \Bbbk(\mathcal E)$, such that $\operatorname{div}(\alpha)=nA-nB$. Let $\mathcal E$ be defined by the equation $y^2= 1 + ax + bx ^ 2 + cx ^ 3 $, such that $x(A)=0$, $y(A)=1$ and $B$ is the infinite point.

\begin{theorem}
Let rational function $\frac{P(x)}{Q(x)}$ be the Pad\'e approximation of $\sqrt{1 + ax + bx ^ 2 + cx ^ 3}$ of order $[\lfloor n/2\rfloor, \lfloor(n-3)/2\rfloor]$. Then for the function $\mathfrak{a}=p(x)-q(x)y$ on $\mathcal E$ we have $\operatorname{div}(\mathfrak a)=(n-1)A+C-B$, where $C$ is some point on $\mathcal E$.
\end{theorem}

\textbf{Proof:} See ~\cite{Grag}. $\Box$

Pad\'e approximation of the function $ f (x) $ is defined by the Taylor coefficients of $ f (x) $. The condition $A=C$ gives us the equation on $a$, $b$, $c$.

\subsection{The case $n=6$}

After series of calculation we get equation for the family of primitive Abel pairs of degree $6$ parameterized by a variable $t$ (i.e. $X_1(6)$ is rational): $${\alpha}^{2}+ \left( -t \left( t-1 \right) ^{2}- \left( 3\,t+1 \right)
 \left( t-1 \right) x-4\,{x}^{2}+2\,{x}^{3} \right) \alpha+{x}^{6}=0$$

$$\varkappa_6(\mathbb C)=-{\frac {1}{2^{14}\cdot 3^{12}}}\,{\frac { \left( 9\,t-1 \right)  \left( 81\,{
t}^{3}-27\,{t}^{2}+99\,t-25 \right) ^{3} \left( 9\,t-25 \right) ^{2}}{{t}^{5} \left( t-1 \right) ^{4}}}$$

Let we see what zeros and poles of $\varkappa_6$ are.

The poles of $\varkappa_6$ are located in the point $t=0$, $t=1$ and $t=\infty$. for $t=0$ equation takes the form $\alpha^2-xR_4(x)\alpha+x^6=0$, At $t=1$ we have $\alpha^2+4x^2R_2(x/4)\alpha+x^6=0$, and in $t=\infty$ we have $\alpha^2+2x^3\alpha+x^6=0$.

The zeros of $\varkappa_6$ are located in $t=1/9$, $t=25/9$ in the roots of the equation $81\,{t}^{3}-27\,{t}^{2}+99\,t-25=0$. Over $t=1/9$ we get the equation $\alpha^2+\frac{2^6}{3^6}R_6(\frac94 x)\alpha+x^6=0$.

Over $t=25/9$ we get Abel-Belyi pair $\Box_{1,1,2,2}$, and over the roots of $81\,{t}^{3}-27\,{t}^{2}+99\,t-25=0$ the Abel-Belyi pairs $\hexagon_{4,1,1}$, $\hexagon_{3,2,1}$, $\hexagon_{3,1,2}$ lie.

We consider the reduction of this family of curves to the characteristic $5$. Let's look at the zeros and poles of $\varkappa_6$. Zero of $\varkappa_6$ from $t=25/9$ now moves to $t=0$, as well as one of the roots of equation $81\,{t}^{3}-27\,{t}^{2}+99\,t-25=0$. In the numerator of $\varkappa_6$ we have $t^5\cdot (t^2+3t+4)^3\cdot(t+1)$. So the pole at $t=0$ disappears:$$\varkappa_6(\overline{\mathbb F_5})=-\frac{(t^2+3t+4)^3\cdot(t+1)}{(t-1)^4}$$

From this calculation we see that over $\overline{\mathbb F_5}$ there are only $2$ primitive Abel-Belyi pairs with the set of valencies $(n|3,1,\ldots,1)$, unlike over $\mathbb C$, where exist $3$ such pairs. Also over $\overline{\mathbb F_5}$ there are no primitive Abel-Belyi pairs, whose set of valencies $(n|2,2,1\ldots,1)$, unlike over $\mathbb C$ where one such pair exist.

So $m_{\hexagon}(6,\mathbb C)-m_{\hexagon}(6,\overline{\mathbb F_5})=1$ and $m_{\Box}(6,\mathbb C)-m_{\Box}(6,\overline{\mathbb F_5})=1$, in accordance with theorem 6.3.

 \end{document}